\documentclass[11pt]{article}
\usepackage[utf8]{inputenc}
\usepackage{amsmath,amssymb,epsfig,bbm}
\usepackage{stmaryrd,mathabx}
\usepackage{comment}
\usepackage{color}
\usepackage[T1]{fontenc}

\usepackage[textsize=small]{todonotes}
\usepackage{enumitem}
\usepackage{varwidth}
\setlist{nolistsep}
\usepackage{hyperref}

\newtheorem{defi}{Definition}[section]
\newtheorem{prop}[defi]{Proposition}
\newtheorem{theo}[defi]{Theorem}
\newtheorem{conj}[defi]{Conjecture}
\newtheorem{lemm}[defi]{Lemma}
\newtheorem{coro}[defi]{Corollary}
\newtheorem{rema}[defi]{Remark}
\newtheorem{exem}[defi]{Example}
\newtheorem{exems}[defi]{Examples}

\newcommand{\bdefi}{\begin{defi}}
\newcommand{\edefi}{\end{defi}}
\newcommand{\bprop}{\begin{prop}}
\newcommand{\eprop}{\end{prop}}
\newcommand{\btheo}{\begin{theo}}
\newcommand{\etheo}{\end{theo}}
\newcommand{\blemm}{\begin{lemm}}
\newcommand{\brema}{\begin{rema}}
\newcommand{\erema}{\end{rema}}
\newcommand{\bexer}{\begin{exem}}
\newcommand{\eexer}{\end{exem}}
\newcommand{\bexems}{\begin{exems}}
\newcommand{\eexems}{\end{exems}}
\newcommand{\bconj}{\begin{conj}}
\newcommand{\econj}{\end{conj}}
\newcommand{\elemm}{\end{lemm}}
\newcommand{\bcoro}{\begin{coro}}
\newcommand{\ecoro}{\end{coro}}
\newcommand{\dem}{\noindent{\bf Proof. }}

\usepackage{mathrsfs}
\renewcommand\mathcal{\mathscr}

\newcommand{\B}{{\cal B}}

\newcommand{\D}{{\cal D}}

\newcommand{\M}{{\cal M}}

\newcommand{\OOO}{{\cal O}}

\newcommand{\V}{{\cal V}}

\newcommand{\maths}[1]{{\mathbb #1}}  

\newcommand{\FF}{\maths{F}}

\newcommand{\LL}{\maths{L}}
\newcommand{\NN}{\maths{N}}

\newcommand{\PP}{\maths{P}}

\newcommand{\RR}{\maths{R}}
\newcommand{\SSS}{\maths{S}}

\newcommand{\ZZ}{\maths{Z}}

\newcommand{\aaa}{{\mathfrak a}}

\newcommand{\mmm}{{\mathfrak m}}

\newcommand{\ppp}{{\mathfrak p}}

\newcommand{\uuu}{{\mathfrak u}}

\newcommand{\zzz}{{\mathfrak z}}

\newcommand{\weakstar}{\overset{*}\rightharpoonup}
\newcommand{\ra}{\rightarrow}
\newcommand{\bs}{\backslash}

\newcommand{\wh}[1]{{\widehat{#1}}}

\newcommand{\ga}{\gamma}
\newcommand{\Ga}{\Gamma}

\newcommand{\Snu}{\SSS^1_\omega}
\newcommand{\Snusharp}{\SSS^{1,\sharp}_\omega}
\newcommand{\prim}{{\rm prim}}
\newcommand{\norm}[1]{\| #1\|_\omega}
\newcommand{\dirp}[1]{\widecheck{#1}}

\newcommand{\cqfd}{\hfill$\Box$}

\newcommand{\bigO}{\operatorname{O}}
\newcommand{\card}{{\operatorname{Card}}}

\newcommand{\id}{\operatorname{id}}

\newcommand{\SL}{\operatorname{SL}}

\newcounter{fig}

\title{Effective equidistribution of lattice points \\ 
in positive characteristic}
\author{Tal Horesh \and Fr\'ed\'eric Paulin} 
\begin{document}
\bibliographystyle{../alphanum}
\maketitle
\begin{abstract}
Given a place $\omega$ of a global function field $K$ over a finite
field, with associated affine function ring $R_\omega$ and completion
$K_\omega$, the aim of this paper is to give an effective joint
equidistribution result for renormalized primitive lattice points
$(a,b)\in {R_\omega}^2$ in the plane ${K_\omega}^2$, and for renormalized
solutions to the gcd equation $ax+by=1$. The main tools are
techniques of Goronik and Nevo for counting lattice points in
well-rounded families of subsets.  This gives a sharper analog in positive
characteristic of a result of Nevo and the first author for the
equidistribution of the primitive lattice points in $\ZZ^2$.
\end{abstract}

\section{Introduction}
\label{sect:intro}

This paper has two motivations. The first one is the following result
of Dinaburg-Sinai \cite{DinSin90}. Given two coprime positive integers
$a<b$, let $(x_{0},y_{0})$ be a shortest (with respect to the supremum
norm $\|(x,y)\|_\infty= \max\{|x|,|y|\}$) solution to the equation
$|ax+by|=1$ with unknown $(x,y) \in\ZZ^2$.  Dinaburg-Sinai proved that
the quotients of norms
$$
\frac{\|(x_{0},y_{0})\|_\infty}{\|(a,b)\|_\infty}
$$ 
equidistribute in the interval $[0,1]$ as $\|(a,b)\|_\infty$ tends
to $+\infty$.

The second motivation is the well-studied Linnik problem of
equidistribution on the unit sphere $\SSS^{n-1}$ of the directions of
integral vectors in the Euclidean space $\RR^{n}$ for $n\geq 2$. See
for instance \cite{Duke88, Schmidt98, Duke03, Duke07, EinLinMicVen11,
  BenOh12, EllMicVen13, AkaEinSha16a, AkaEinSha16b} as well as the
joint works of the first author \cite{HorNev16,HorKar19}. Let us
denote by $\ZZ^{n}_{\rm prim}$ the set of primitive integral vectors,
by $\operatorname{Leb}_{\SSS^{n-1}}$ the spherical measure on
$\SSS^{n-1}$ renormalized to be a probability measure, and by
$\Delta_x$ the unit Dirac mass at any point $x$. A simple version of
this equidistribution phenomenon is the now well-known fact that, as
$N\ra+\infty$,
$$
\frac{1}{\card\{v\in\ZZ^{n}_{\rm prim}: \|v\| \leq N\}}\;\;
\sum_{v\in\ZZ^{n}_{\rm prim}\;:\;\|v\| \leq N} \;\Delta_{\frac{v}{\|v\|}}
\;\;\weakstar\;\;\operatorname{Leb}_{\SSS^{n-1}}
$$ 
for the weak-star convergence of measures on the compact space
$\SSS^{n-1}$.  A connection between the two motivations is that when $n=2$,
an integral vector $(a,b)$ is primitive if and only if there exists an
integral vector $(x,y)$ with $|ax+by|=1$.

\medskip The goal of this paper is to address analogous questions in
local fields with positive characteristic. In this introduction, we
describe our results in the special following case.

Let $\FF_q$ be a finite field of order a positive power $q$ of some
positive prime, and let $K=\FF_q(Y)$ be the field of rational
functions in one variable $Y$ over $\FF_q$. Let $R=\FF_q[Y]$ be the
ring of polynomials in $Y$ over $\FF_q$, let $\wh{K}=\FF_q((Y^{-1}))$
be the non-Archimedean local field of formal Laurent series in
$Y^{-1}$ over $\FF_q$ and let $\OOO=\FF_q[[Y^{-1}]]$ be the local ring of
$\wh{K}$ (consisting of formal power series in $Y^{-1}$ over
$\FF_q$). We denote by $|\cdot|$ the complete nonarchimedean absolute
value on $\wh{K}$ such that $|P|=q^{\deg P}$ for every $P\in R$.

We endow $\wh{K}$ with its Haar measure $\mu_{\wh{K}}$ standardly
normalized so that $\mu_{\wh{K}}(\OOO)=1$, and the quotient $\wh{K}/R$
with the induced measure $\mu_{\wh{K}/R}$. We also endow the plane
$\wh{K}^{\,2}$ with the product measure and with the supremum norm. We
denote by $\SSS^1_\infty$ the (compact-open) unit sphere of
$\wh{K}^{\,2}$, that we equip with the restriction
$\mu_{\SSS^1_\infty}$ of the product measure.

Given $v=(a,b)\in \wh{K}^{\,2}-\{0\}$, we denote by $\|v\|_\infty=
\max\{|a|,|b|\}\in q^\ZZ$ its supremum norm, and by $z_v=a$ if
$|a|\geq|b|$, and $z_v=b$ otherwise, its component with maximum
absolute value. We also denote by $\dirp v=\big(a\,
Y^{-\log_q\|v\|_\infty}, b\, Y^{-\log_q\|v\|_\infty}\big)$ the vector
$v$ canonically renormalised to be in the unit sphere $\SSS^1_\infty$,
of which we think as the {\it direction} of $v$.

We let $R^2_\prim$ denote the set of elements $v=(a,b)$ in the
standard $R$-lattice $R^{\,2}$ of the plane $\wh{K}^{\,2}$ that are
{\it primitive}, that is, satisfy $aR+bR=R$. Let $w_v=(-y',x')$ be such
that $(x',y')$ is a solution to the gcd equation $ax+by=1$ of $v$,
with unknown $(x,y)\in R^{\,2}$. We could for instance take the {\it
  shortest} one, that is, the one with the smallest supremum norm (see
Section \ref{sect:cfe} for the existence and uniqueness). What follows
is actually independent of the choice of $w_v$.

\medskip
The following result is a joint equidistribution theorem, with error
term, for the direction and renormalized gcd solution of
the primitive lattice points in the nonarchimedean plane
$\wh{K}^{\,2}$.

Error terms in equidistribution results usually require smoothness
properties on test functions. The appropriate smoothness regularity of
functions defined on totally disconnected spaces like $\wh K^N$ for
$N\in\NN$ is the locally constant one. For every metric space $E$ and
$\epsilon > 0$, a bounded map $f : E \ra \RR$ is {\it
  $\epsilon$-locally constant} if it is constant on every closed ball
of radius $\epsilon$ in $E$. Its $\epsilon$-locally constant norm is
$\|f\|_\epsilon= \frac{1}{\epsilon}\sup_{x\in E}|f(x)|$.

\btheo\label{theo:mainintro} For the weak-star
convergence on the compact space $\SSS^1_\infty\times (\wh K/R)$, we
have, as $n\ra+\infty$,
$$
\frac{1}{q^2(q-1)} \;q^{-2n}
\sum_{v\in R^2_{\prim}\;:\;\|v\|_\infty= q^n} 
\Delta_{\dirp v}\otimes\Delta_{\frac{z_{w_{v}}}{z_v}+R}\;\;\weakstar\;\;
\mu_{\SSS^1_\infty}\otimes \mu_{\wh K/R}\;.
$$ 
Furthermore, there exists $\tau\in\;]0,\frac{1}{8}]$ such that for
all $\epsilon,\delta>0$, there is a mutiplicative error term of the form
$1+\bigO_{\delta}(q^{2n(-\tau+\delta)}\,\|f\|_\epsilon\,\|g\|_\epsilon)$
when evaluated on $(f,g)$ for all $\epsilon$-locally constant maps
$f:\SSS^1_\infty\ra\RR$ and $g:\wh K/R\ra\RR$.
%
%
\etheo

The constant $\tau$ is explicit in terms of representation-theoretic
data for the locally compact group $\SL_2(\wh K)$. We will actually
prove a more general version of this result, when $\wh K$ is replaced
by any (global) function field in one variable over a finite field and
when congruence properties are added, see Theorem \ref{theo:main}.
See also Corollary \ref{coro:counting} for a counting corollary of
primitive lattice points.

We begin in Subsection \ref{subsect:functionfields} by recalling basic
facts about functions fields over finite fields. In Subsection
\ref{subsect:modulargroup}, we define the various closed subgroups of
the totally disconnected locally compact group $\SL_2(\wh K)$ which
will be useful in order to transfer arithmetic information on lattice
points in the plane to group-theoretic information. We will also
discuss the properties of their Haar measures. In Section
\ref{sect:correspondance}, we give a precise correspondence between
primitive lattice points and elements in the Nagao-Weyl modular group
$\SL_2(\FF_q[Y])$. We adapt in Section \ref{sect:equidistrib} the
results of Gorodnik-Nevo \cite{GorNev12} on counting lattice points in
well-rounded subsets of semi-simple Lie groups, and check that a
family of nice compact-open subsets coming from a mixture of the LU
and Iwasawa decompositions of $\SL_2(\wh K)$ is indeed well-rounded.
Finally, in Section \ref{sect:cfe}, we give an application to the
distribution properties of the continued fraction expansions of
elements in $\FF_q(Y)$, thus giving an analogue to the result of
Dinaburg-Sinai in \cite{DinSin90} described in the beginning of this
Introduction.

\bigskip
{\small \noindent {\it Acknowledgments. } The authors warmly thank
  Amos Nevo for having presented one to the other in a beautiful
  conference in Goa in February 2016, where the idea of this paper was
  born. The first author thanks the IHES for two post-doctoral years
  when most of this paper was discussed, and the Topology team in
  Orsay for financial support at the final stage.}

\section{Background on function fields and their modular groups}
\label{sect:background}

\subsection{Global function fields}
\label{subsect:functionfields}

We refer for instance to \cite{Goss96,Rosen02} and
\cite[Chap.~14]{BroParPau19} for the content of this Section.

Let $\FF_q$ be a finite field of order $q$, where $q$ is a positive
power of a positive prime. Let $K$ be a (global) {\it function field}
over $\FF_q$, that is, the function field of a geometrically connected
smooth projective curve ${\bf C}$ over $\FF_q$, or equivalently an
extension of $\FF_q$ of transcendance degree $1$, in which $\FF_q$ is
algebraically closed. We denote by ${\bf g}$ the genus of the curve ${\bf
  C}$.

There is a bijection between the set of closed points of ${\bf C}$ and
the set of (normalised discrete) valuations $\omega$ of its function
field $K$, where the valuation of a given element $f\in K$ is the order
of the zero or the opposite of the order of the pole of $f$ at the
given closed point. We fix such an element $\omega$ from now on.

We denote by $K_\omega$ the completion of $K$ for the valuation
$\omega$, and by 
$$
\OOO_\omega= \{x\in K_\omega\;:\;\omega(x)\geq 0\}
$$
the valuation ring of (the unique extension to $K_\omega$) of
$\omega$.  Let us fix a uniformiser $\pi_\omega\in K_\omega$ of
$\omega$, that is, an element in $K_\omega$ with
$\omega(\pi_\omega)=1$.  We denote by $q_\omega$ the order of the
residual field $\OOO_\omega/\pi_\omega\OOO_\omega$ of $\omega$, which
is a (possibly proper) power of $q$.  We normalize the absolute value
associated with $\omega\,$ as usual: for every $x\in K_\omega$, we
have the equality
$$
|\,x\,|_\omega=(q_\omega)^{-\omega(x)}\;.
$$

Finally, let $R_\omega$ denote the affine algebra of the affine curve
${\bf C} -\{\omega\}$, consisting of the elements of $K$ whose only
poles are at the closed point $\omega$ of ${\bf C}$. Its field of
fractions is equal to $K$.

\medskip
The case in the introduction corresponds to ${\bf C}=\PP^1$ and
$\omega=\omega_\infty$ the valuation associated with the point at infinity
$[1:0]$. Then

$\bullet$~ $K=\FF_q(Y)$ is the field of rational functions in one
variable $Y$ over $\FF_q$, 

$\bullet$~ $\omega_\infty$ is the valuation defined, for all $P,Q\in
\FF_q[Y]$, by
$$
\omega_\infty(P/Q)=\deg Q-\deg P\;.
$$ 

$\bullet$~ $R_{\omega_\infty}=\FF_q[Y]$ is the (principal) ring of
polynomials in one variable $Y$ over $\FF_q$,

$\bullet$~ $K_{\omega_\infty}= \FF_q((Y^{-1}))$ is the field of formal
Laurent series in one variable $Y^{-1}$ over $\FF_q$,

$\bullet$~ $\OOO_{\omega_\infty}= \FF_q[[Y^{-1}]]$ is the ring of formal
power series in one variable $Y^{-1}$ over $\FF_q$, $\pi_{\omega_\infty}
=Y^{-1}$ is the usual choice of an uniformizer, and $q_{\omega_\infty}=q$.

\medskip
Recall (see for instance \cite[II.2 Notation]{Serre83}) that $R_\omega$
is a Dedeking ring, not principal in general.  By for instance
\cite[Eq.~(14.2)]{BroParPau19}, we have
\begin{equation}\label{eq:interRnupiOnu} 
R_\omega \cap \OOO_\omega=\FF_q\;.
\end{equation}

\blemm\label{lem:maxmin} 
For all $a,b,c,d\in K_\omega-\FF_q$ such that $ad-bc=1$ and
$|\,a\,|_\omega\geq|\,b\,|_\omega$, we have $|\,c\,|_\omega\geq|\,d\,|_\omega$.
\elemm

\dem 
The equality $ad-bc=1$ implies that $\omega(ad-bc)=0$. We have
$\omega(ad)<0$ and $\omega(bc)<0$ since the only elements of $R_\omega$ which
have nonnegative valuations are the elements in the ground field
$\FF_q$ by Equation \eqref{eq:interRnupiOnu}. Therefore
$\omega(ad)=\omega(bc)$ and
$$
\omega(c)-\omega(d)=\omega(a)-\omega(b)\;.
$$
The left-hand side is nonpositive, since the right-hand side
is. This proves the result. 
\cqfd

\medskip
The (absolute) {\it norm} of a nonzero ideal $I$ of the ring $R_\omega$ is
$N(I)=[R_\omega:I]=|R_\omega/I|$. {\em Dedekind's zeta function} of $K$ is (see
for instance \cite[\S 7.8]{Goss96} or \cite[\S 5]{Rosen02})
$$
\zeta_K(s)=\sum_I \frac{1}{N(I)^s}
$$ 
where the summation runs over the nonzero ideals $I$ of $R_\omega$. By
(for instance) \cite[\S 5]{Rosen02}, it is a rational function of
$q^{-s}$ with simple poles at $s=0,s=1$. In particular, when
$K=\FF_q(Y)$, then (see \cite[Theo.~5.9]{Rosen02})
\begin{equation}\label{zetamoinsun}
\zeta_{\FF_q(Y)}(-1)=\frac{1}{(q-1)(q^2-1)}\;.
\end{equation}

We denote by 
$$
R_{\omega,\prim}^2=\{(a,b)\in{R_{\omega}}^2:aR_\omega+bR_\omega=R_\omega\}
$$ 
the set of {\it primitive} elements in the lattice ${R_{\omega}}^2$ in the
plane ${K_\omega}^2$.  Note that since $R_\omega$ is not always principal, not
every point of ${R_{\omega}}^2$ is an $R_\omega$-multiple of an element of
$R_{\omega,\prim}^2$.

\medskip
For every $v\in {K_\omega}^2-\{0\}$, we write $v=(x_v,y_v)$, and define
\begin{equation}\label{eq:znu}
z_v=\left\{\begin{array}{ll }x_v& {\rm if}\;|\,x_v|_\omega\geq|\,y_v|_\omega\\
y_v& {\rm if}\;|\,x_v|_\omega<|\,y_v|_\omega\end{array}\right.
\;\;\;{\rm and}\;\;\;
z'_v=\left\{\begin{array}{ll }y_v& {\rm if}\;|\,x_v|_\omega\geq|\,y_v|_\omega\\
x_v& {\rm if}\;|\,x_v|_\omega<|\,y_v|_\omega\;,\end{array}\right.
\end{equation} 
%
%
as well as
\begin{equation}\label{eq:dirpv}
\|v\|_\omega=\max\{|\,x_v|_\omega,|\,y_v|_\omega\},\;\;\;
v^\perp=(y_v,-x_v)\;\;\;{\rm and}\;\;\;\dirp v =
\pi_\omega^{\log_{q_\omega}(\norm{v})}\,v\;.
\end{equation} 
We denote the unit sphere in the plane ${K_\omega}^2$ endowed with the
supremum norm $\|\cdot\|_\omega$ by
$$
\Snu=\{v\in{K_\omega}^2: \|v\|_\omega=1\}\;.
$$

Note that $v^\perp$ has the same norm as $v$ and belongs to
$R_{\omega,\prim}^2$ if $v$ does, and that $\Snu= \{\,{\dirp v}:
v\in{K_\omega}^2-\{0\}\,\}$.  We think of $\dirp v$ as the {\it
  direction} (or renormalisation) of $v$, it is a prefered element in
the intersection of the unit sphere $\Snu$ with the vector line
defined by $v$.

\medskip
We denote by $\|\mu\|$ the total mass of any finite measure $\mu$.  We
let $\mu_{K_\omega}$ denote the Haar measure of the (abelian) locally
compact topological group $(K_\omega,+)$, normalised so that
$\mu_{K_\omega}(\OOO_\omega) =1$.  This measure scales as follows
under multiplication: for all $\lambda,x\in K_\omega$, we have
\begin{equation}\label{eq:homothetyhaar}
d\mu_{K_\omega}(\lambda x)=|\lambda|_\omega\;d\mu_{K_\omega}(x)\;.
\end{equation} 
We denote by $\mu_{K_\omega/R_\omega}$ the induced Haar measure on the
compact additive topological group $K_\omega/R_\omega$. Using the above
scaling for the first equation and \cite[Lem.~14.4]{BroParPau19} for the
second one,  for every $m\in\NN$, we have the equality
\begin{equation}\label{eq:haarKmodR}
\mu_{K_\omega}(\pi_\omega^m\OOO_\omega)=q_\omega^{-m}\;\;\;{\rm and}\;\;\;
\|\,\mu_{K_\omega/R_\omega}\|=q^{{\bf g}-1}\;.
\end{equation}

We endow ${K_\omega}^2$ with the product $\mu_{K_\omega}\otimes
\mu_{K_\omega}$ of the Haar measures on each factor. Note that the
unit ball of ${K_\omega}^2$ is ${\OOO_\omega}^2$, so that for every
$k\in\ZZ$, the measure of any ball in ${K_\omega}^2$ of radius
$q_\omega^k$, which is of the form $v+\pi_\omega^{-k}{\OOO_\omega}^2$
for some $v\in {K_\omega}^2$, is equal to $q_\omega^{2k}$.

We denote by $\mu_{\Snu}$ the restriction to the compact-open subset
$\Snu$ of ${K_\omega}^2$ of the product measure. Since
\begin{equation}\label{eq:mesoootimes}
\mu_{K_\omega}(\OOO_\omega^\times)=
\mu_{K_\omega}(\OOO_\omega -\pi_\omega\OOO_\omega)=1-q_\omega^{-1}
\end{equation}
by Equation \eqref{eq:haarKmodR}, and since $\Snu=(\OOO_\omega^\times
\times \OOO_\omega)\cup (\OOO_\omega\times\OOO_\omega^\times)$, the
total mass of $\mu_{\Snu}$ is
\begin{equation}\label{eq:mesSnu}
  \|\mu_{\Snu}\|= (1-q_\omega^{-1})+(1-q_\omega^{-1})-(1-q_\omega^{-1})^{2}=
  \frac{q_\omega^{2}-1}{q_\omega^2}\;.
\end{equation}

\subsection{The modular group}
\label{subsect:modulargroup}

The aim of this section is to introduce the various closed subgroups
of the special linear group of the plane ${K_\omega}^2$ that will be
useful in order to transfer arithmetic information concerning lattice
points in ${R_\omega}^2$ into group-theoretic information. We will also
discuss the properties of their Haar measures.

Let $G=\SL_2(K_\omega)$, which is a totally disconnected locally
compact topological group. The {\em modular group} $\Ga=
\SL_2(R_\omega)$ is a non-uniform lattice in $G$. When ${\bf C}=
\PP^1$ and $\omega=\omega_\infty$ as in the Introduction, then up to finite
index, it is called {\it Nagao's lattice} (see \cite{Nagao59,Weil70}).
For every nonzero ideal $I$ of $R_\omega$, we denote by $\Ga_0[I]$ the
{\it Hecke congruence subgroup} of $\Ga$ modulo $I$:
$$
\Ga_0[I]=\big\{\big(\begin{smallmatrix} a & c \\ b & d
\end{smallmatrix}\big) \in \Ga :b\in I\big\}\;.
$$
By \cite[Lem.~16.5]{BroParPau19}, the index of $\Ga_0[I]$ in $\Ga$ is
\begin{equation}\label{eq:indexHecke}
  \big[\;\Ga:\Ga_0[I]\,\big]=
  N(I)\;\prod_{\ppp | I} \Big(1+\frac{1}{N(\ppp)}\Big)\;.
\end{equation}
where the product ranges over the prime factors $\ppp$ of the ideal $I$.

For every closed subgroup $H$ of $G$, we denote by $H(\OOO_\omega)$ the
compact-open subgroup $H\cap \M_2(\OOO_\omega)$ of $H$, and by $\mu_H$ the
(left) Haar measure of $H$ normalized so that
$$
\mu_H(H(\OOO_\omega))=1\;.
$$

Note that $G$ is unimodular. For every lattice $\Ga'$ of $G$, we
denote by $\mu_{\Ga'\bs G}$ the measure on $\Ga'\bs G$ induced by
$\mu_G$. By Exercice 2 e) in \cite[II.2.3]{Serre83} (which normalizes
the Haar measure of $G$ so that the mass of $G(\OOO_\omega)$ is
$q_\omega-1$), the total mass of $\mu_{\Ga\bs G}$ is
\begin{equation}\label{eq:covolGa}
\|\mu_{\Ga\bs G}\|=\zeta_K(-1)\;.
\end{equation}

Let $Z$ be the diagonal subgroup of $G$, let $U^-$ and $U^+$ be its
lower and upper unipotent triangular subgroups, and let $P^-=U^-Z$ be
its lower triangular Borel subgroup. We also consider the Cartan
subgroup $A=\big\{\big(\begin{smallmatrix} \pi_\omega^n & 0 \\ 0 &
  \pi_\omega^{-n}\end{smallmatrix} \big):n\in\ZZ\big\}$ of $G$, whose
centralizer in $G$ is $Z$.

Since $A(\OOO_\omega)=\{\id\}$ has measure one, the Haar measure
$\mu_A$ on $A$ is exactly the counting measure:
\begin{equation}\label{eq:muA}
\mu_A=\sum_{g\in A} \;\Delta_g\;.
\end{equation}

The maps from $K_\omega$ to $U^-$ and $U^+$, defined respectively by
$\alpha\mapsto \uuu^-(\alpha)= \big(\begin{smallmatrix}
1 & 0\\ \alpha & 1\end{smallmatrix} \big)$ and $\alpha\mapsto
\uuu^+(\alpha)= \big(\begin{smallmatrix} 1 & \alpha \\ 0 &
1\end{smallmatrix} \big)$, are homeomorphisms (and even abelian group
isomorphisms). They send $\OOO_\omega$ to $U^\pm(\OOO_\omega)$, and
the Haar measure of $(K_\omega,+)$ to the Haar measure of $U^\pm$: namely, for
(almost) all $\alpha\in K_\omega$, we have
\begin{equation}\label{eq:muUpm}
d\mu_{U^\pm}(\uuu^\pm(\alpha))=d\mu_{K_\omega}(\alpha)\;.
\end{equation}
Similarly, the map from the multiplicative group ${K_\omega}^\times$
to $Z$, defined by $\alpha\mapsto\big(\begin{smallmatrix} \alpha & 0
  \\0 & \alpha^{-1}\end{smallmatrix} \big)$, is a homeomorphism (and
even an abelian group isomorphism). It sends ${\OOO_\omega}^\times$ to
$Z(\OOO_\omega)$, and the restriction to ${K_\omega}^\times$ of the
Haar measure $\mu_{K_\omega}$ to a multiple of the Haar measure of
$Z$: namely, for (almost) all $\alpha\in {K_\omega}^\times$, by
Equation \eqref{eq:mesoootimes}, we have
\begin{equation}\label{eq:muZ}
\frac{q_\omega-1}{q_\omega}\;d\mu_{Z}(\big(\begin{smallmatrix} \alpha & 0
  \\ 0 & \alpha^{-1}\end{smallmatrix} \big))=
d\mu_{K_\omega}(\alpha)\;.
\end{equation}

Let 
$$
\Snusharp=\big\{v\in \Snu\;:\; |\,x_v|_\omega\geq |\,y_v|_\omega\big\} =
\OOO_\omega^\times\times\OOO_\omega\;,
$$ which is a compact-open subset of the plane ${K_\omega}^2$. The map
from $\Snusharp$ to $P^-(\OOO_\omega)$ defined by $(\alpha,\beta)
\mapsto \ppp^-(\alpha,\beta)= \big(\begin{smallmatrix} \alpha & 0
  \\ \beta & \alpha^{-1}\end{smallmatrix} \big)$ is a homeomorphism.
It sends the restriction to $\Snusharp$ of the measure $\mu_{\Snu}$ to
a multiple of the Haar measure of $P^-(\OOO_\omega)$: since
$\mu_{P^-(\OOO_\omega)}$, $\mu_{U^-(\OOO_\omega)}$ and
$\mu_{Z(\OOO_\omega)}$ are probability measure, by Equations
\eqref{eq:muUpm} and \eqref{eq:muZ}, we have, for (almost) every
$\alpha\in \OOO_\omega^\times$ and $\beta\in\OOO_\omega$
\begin{align}\nonumber
  d\mu_{P^-(\OOO_\omega)}(\ppp^-(\alpha,\beta))&=
  d\mu_{U^-(\OOO_\omega)}(\uuu^-(\beta))\;d\mu_{Z(\OOO_\omega)}( 
\big(\begin{smallmatrix} \alpha & 0
  \\ 0 & \alpha^{-1}\end{smallmatrix} \big))=
\frac{q_\omega}{q_\omega-1}\;d\mu_{K_\omega}(\alpha)\,d\mu_{K_\omega}(\beta)
\\ &=\frac{q_\omega}{q_\omega-1}\;d\mu_{\Snu}(\alpha,\beta)\;.
\label{eq:decompmes}
\end{align}

\medskip
We will need the following refined LU decomposition of elements of the
special linear group $G$.  Let $g=\left(\begin{smallmatrix}\alpha &
  \gamma\\ \beta & \delta \end{smallmatrix}\right)\in G$ with $\alpha
\neq 0$. Then there are unique elements $\uuu^\pm_g\in U^\pm$,
$\mmm_g\in Z(\OOO_\omega)$ and $\aaa_g\in A$ such that
$$
g= \uuu^-_g\;\mmm_g\;\aaa_g\;\uuu^+_g\;.
$$
We have 
\begin{equation}\label{eq:refineLU}
\uuu^-_g=
\begin{pmatrix} 1 & 0\\\frac{\beta}{\alpha} & 1\end{pmatrix},\;\;
\uuu^+_g=
\begin{pmatrix} 1 &\frac{\ga}{\alpha}\\ 0 & 1\end{pmatrix},\;\;
\mmm_g=
\begin{pmatrix} \alpha\pi_\omega^{-\omega(\alpha)} & \!\!\!0\\
\!\!\!0 & \!\!\!\!\!\!\alpha^{-1}\pi_\omega^{\omega(\alpha)}
\end{pmatrix},\;\;
\aaa_g=
\begin{pmatrix} \pi_\omega^{\omega(\alpha)} & \!\!0\\\!\!0 & 
\!\!\!\!\pi_\omega^{-\omega(\alpha)}\end{pmatrix}.
\end{equation}
We also consider
\begin{equation}\label{eq:pppg}
\ppp_g=\uuu^-_g\;\mmm_g =
\begin{pmatrix} \alpha\pi_\omega^{-\omega(\alpha)} & 0\\
\beta\pi_\omega^{-\omega(\alpha)} & \alpha^{-1}\pi_\omega^{\omega(\alpha)}
\end{pmatrix}\in P^-\;.\end{equation}
Note that if $\omega(\alpha)\leq \omega(\beta)$, that is, if
$|\,\alpha\,|_\omega\geq |\,\beta\,|_\omega$, then $\ppp_g\in
P^-(\OOO_\omega)=U^-(\OOO_\omega)Z(\OOO_\omega)$, so that $\ppp_g$
belongs to the maximal compact subgroup $G(\OOO_\omega)$ of $G$. In
particular, the writing $g= \ppp_g\; \aaa_g\;\uuu^+_g$ is an Iwasawa
decomposition of $g$.

We conclude this section by providing the expression for the Haar measure
of $G$ in the refined LU decomposition. The product map from
$U^-\times Z(\OOO_\omega)\times A\times U^+$ to $G$ is an homeomorphism
onto an open-dense subset with full Haar measure in $G$, and the
following result says that the Haar measure of $G$ is absolutely
continuous with respect to the product of the Haar measures of the
factors. The main point of its proof is to compute the constant. We
denote by $\chi:Z\ra {K_\omega}^\times$ the standard character
$\big(\begin{smallmatrix} \alpha^{-1} & 0\\ 0 & \alpha
\end{smallmatrix} \big)\mapsto \alpha$. 
%
%
It is well known (by the standard action of a split torus on its root
groups) that for all $\zzz\in Z$ and $\alpha\in K_\omega$, we have
\begin{equation}\label{eq:dilatlarate}
\zzz\,\uuu^-(\alpha)\,\zzz^{-1}=\uuu^-(\,\chi(\zzz)^2\,\alpha)
\;\;\;{\rm and}\;\;\;
\zzz^{-1}\,\uuu^+(\alpha)\,\zzz=\uuu^+(\,\chi(\zzz)^2\,\alpha)\,.
\end{equation}

\blemm \label{lem:decompHaarG}
For $\mu_G$-almost every $g\in G$, we have
$$
d\mu_G(g)= \frac{q_\omega}{q_\omega+1}\;|\chi(\aaa_g)|_\omega^{\,-2}\,
d\mu_{U^-}(\uuu^-_g)\;d\mu_{Z(\OOO_\omega)}(\mmm_g)\;
d\mu_{A}(\aaa_g)\;d\mu_{U^+}(\uuu^+_g)\;.
$$
\elemm

\dem By \cite[\S III.1]{Lang75}, since $G$ and $U^+$ are unimodular,
there exists a constant $c_1>0$ such that $d\mu_G(p^-u^+) = c_1\;
d\mu_{P^-}(p^-)\;d\mu_{U^+}(u^+)$ for (almost) every $p^-\in P^-$ and
$u^+\in U^+$, using the product map $P^-\times U^+\ra G$. Note that
$U^-$ is unimodular and that $Z$ normalizes $U^-$ as made precise in
Equation \eqref{eq:dilatlarate}.
Hence there exists a constant $c_2>0$ such that, for (almost) every
$u^-\in U^-$ and $z\in Z$,
$$
|\chi(z)|_\omega^{\,-2}\; d\mu_{U^-}(u^-)\; d\mu_{Z}(z)=c_2\; d\mu_{P^-}(u^-z)
$$ 
This indeed follows by uniqueness from the fact that the left-hand
side defines a Haar measure on $P^-$ using the product map
$(u^-,z)\mapsto u^-z$ from $U^-\times Z$ to $P^-$ (which is an
homeomorphism), by Equations \eqref{eq:homothetyhaar} and
\eqref{eq:muUpm}. Since $Z=Z(\OOO_\omega)A$ with $A$ and $Z(\OOO_\omega)$
abelian and commuting, this proves that there exists a constant
$c_3>0$ such that
\begin{equation}\label{eqctrois}
d\mu_G(g)= c_3\;|\chi(\aaa_g)|_\omega^{\,-2}\;
d\mu_{U^-}(\uuu^-_g)\;d\mu_{Z(\OOO_\omega)}(\mmm_g)\;
d\mu_{A}(\aaa_g)\;d\mu_{U^+}(\uuu^+_g)\;.
\end{equation} 

In order to compute the constant $c_3$, we evaluate the measures on
both sides on the compact-open subgroup $H=\{\big(\begin{smallmatrix}
  \alpha & \ga\\ \beta & \delta \end{smallmatrix} \big)\in
G(\OOO_\omega): \alpha,\delta\in 1+\pi_\omega\OOO_\omega,\; \beta,\ga
\in \pi_\omega\OOO_\omega\}$. This group, being the kernel of the
reduction modulo $\pi_\omega\OOO_\omega$, has index
$|\SL_2(\FF_{q_\omega})| = q_\omega(q_\omega^2-1)$ in
$G(\OOO_\omega)$.  Since $\mu_G(G(\OOO_\omega))=1$, the group $H$ has
Haar measure $\mu_G(H)=\frac{1}{q_\omega(q_\omega^2-1)}$.  By Equation
\eqref{eq:refineLU}, the refined LU decomposition identifies $H$ with
the product $H_{U^-}\times H_{Z}\times H_{U^+}$ in $U^-\times Z\times
U^+$, where
$$
H_{U^-}=
\{\big(\begin{smallmatrix} 1 & 0\\ \beta & 1 \end{smallmatrix} 
\big): \beta\in\pi_\omega\OOO_\omega\},\;
H_{Z}= \{\big(\begin{smallmatrix} \alpha & 0\\ 0 & \alpha^{-1} 
\end{smallmatrix} \big): \alpha\in 1+\pi_\omega\OOO_\omega\},\;
H_{U^+}=\{\big(\begin{smallmatrix} 1 & \ga\\ 0 & 1 \end{smallmatrix} 
\big): \ga\in\pi_\omega\OOO_\omega\}\,.
$$
These groups have index respectively $q_\omega$, $\big|\OOO_\omega^\times/
(1+\pi_\omega\OOO_\omega)\big|=|\FF_{q_\omega}^\times|=q_\omega-1$ and $q_\omega$ in
$U^-(\OOO_\omega)$, $Z(\OOO_\omega)$ and $U^+(\OOO_\omega)$. Hence the measure
of $H$ for the measure on the right-hand side of Equation
\eqref{eqctrois} is equal to $\frac{c_3}{q_\omega^2(q_\omega-1)}$. This
implies that $c_3= \frac{q_\omega}{q_\omega+1}$, as wanted.  
\cqfd

\section{Primitive lattice points seen in the modular group}
\label{sect:correspondance}

Let $K$ be a function field over $\FF_q$, let $\omega$ be a (normalized
discrete) valuation of $K$, let $K_\omega$ be the associated completion
of $K$, and let $R_\omega$ be the affine function ring associated with
$\omega$ (see Section \ref{subsect:functionfields}). The aim of this
section is to naturally associate elements in the modular group
$\Ga=\SL_2(R_\omega)$ to primitive lattice points in ${R_\omega}^2$.

We start by introducing subsets of the plane ${K_\omega}^2$ and of the
group $G=\SL_2(K_\omega)$ which will be technically useful. Let
$$
G^\sharp=\big\{\big(\begin{smallmatrix}\alpha & \ga\\ \beta & \delta
\end{smallmatrix}\big)
\in G\;:\; |\,\alpha\,|_\omega\geq |\,\beta\,|_\omega\big\}
\;\;\;{\rm and}\;\;\;\Ga^\sharp=\Ga\cap G^\sharp \,,$$
$$
K_\omega^{2,\sharp}=\big\{(a,b)\in {K_\omega}^2\;:\; 
|\,a\,|_\omega\geq |\,b\,|_\omega\big\} \;\;\;{\rm and}\;\;\;
R_{\omega,\prim}^{2,\sharp}=R_{\omega,\prim}^{2}\cap K_\omega^{2,\sharp}\;.
$$ 
We identify any element of ${K_\omega}^2$ with the column matrix of
its components. For all measurable subsets $\Theta$ of $\Snu$ and
$\D'$ of $K_\omega$, and for every $n\in\ZZ$, let
\begin{align*}
P^-_\Theta&=\big\{\begin{pmatrix}v' & w'
\end{pmatrix}\in P^-(\OOO_\omega):v'\in\Theta\big\},\\
A_{n}&=\big\{
\big(\begin{smallmatrix}\pi_\omega^{-n} & 0\\ 0 & \pi_\omega^{n}
\end{smallmatrix}\big)\big\}\subset A,\\
U^+_{\D'}&=\big\{ \big(\begin{smallmatrix}1 & \ga\\ 0 & 1
\end{smallmatrix}\big)\in U^+:\ga\in\D'\big\}\,.
\end{align*}
By Lemma \ref{lem:decompHaarG} and the various explicitations of
Haar measures in Equations \eqref{eq:decompmes}, \eqref{eq:muA} and
\eqref{eq:muUpm}, we have
\begin{align}
  \mu_G(P^-_\Theta A_{n}U^+_{\D'}) &=
  \frac{q_\omega}{q_\omega+1} \,\frac{q_\omega}{q_\omega-1}\,
\mu_{\Snu}(\Theta)\,\big({|\pi_\omega^n|_\omega}^{-2}\big)\,
\mu_{K_\omega}(\D')\nonumber \\ &=
\frac{q_\omega^{2n+2}}{q_\omega^2-1}\;\mu_{\Snu}(\Theta)\;\mu_{K_\omega}(\D')\;.
\label{eq:calchaarPThetaA0nUDom}
\end{align}

\medskip
The following result gives a precise 1-to-1 correspondence between
primitive lattice points in $R_{\omega,\prim}^{2,\sharp}$ and
appropriate matrices in the modular group $\Ga$.

\bprop\label{prop:matrepprimvect} Let $\D$ be a fixed (strict)
fundamental domain for the lattice $R_\omega$ acting by translations
on $K_\omega$.  There exists a bijection from
$R_{\omega,\prim}^{2,\sharp}$ to $\Ga^\sharp\cap(P^-\,U^+_\D)$ of the
form $v\mapsto\ga_{v}= \begin{pmatrix} v & w_{v}\end{pmatrix}$ such
that for every $n$ in $\ZZ$, for all measurable subsets $\Theta$
of $\Snu$ and $\D^{\prime}$ of $\D$, and for every nonzero ideal $I$
of $R_\omega$, the following two assertions are equivalent:

\smallskip
\begin{enumerate}
\item[(1)] \label{enu:bijection1} the lattice point $v$ satisfies
  $\norm v= q_\omega^{n}$, $y_v\in I$, $\dirp{v}\in\Theta$ and
  $\frac{x_{w_{v}}}{x_{v}}\in\D'$,
\item[(2)] \label{enu:bijection2} the modular matrix $\ga_{v}$
  belongs to the Hecke congruence subgroup $\Ga_0[I]$ and satisfies
  $\ga_{v}\in P^-_{\Theta}\;A_{n}\;U^+_{\D'}$.
\end{enumerate}
\eprop

\dem Let $v=(a,b)\in R_{\omega,\prim}^{2,\sharp}$. In particular $a\neq
0$ and $\norm{v}=|\,a\,|_\omega$.

In this proof, by {\it solution}, we mean a solution to the equation
$ax+by=1$ with unknown $(x,y)$ varying in ${R_\omega}^2$. Given a {\it
  solution} $w_0=(x_0,y_0)$, we claim that the set of {\it solutions}
is $\{w_0+ \lambda\, v^\perp\;:\; \lambda\in R_\omega\}$, where $w\mapsto
w^\perp$ is defined in Section \ref{subsect:functionfields}. Indeed,
any other {\it solution} $(x,y)\neq (x_0,y_0)$ satisfies $a(x-x_0)=
b(y_0-y)$. We may assume that $b\neq 0$, since otherwise $a\in
R_\omega^\times$ and the result is clear. Then $x\neq x_0$ and $y\neq
y_0$, so that the nonzero principal ideal $(a)$, being coprime with
the principal ideal $(b)$ in the Dedeking ring $R_\omega$, divides the
principal ideal generated by $y_0-y$, and $y-y_0$ is a multiple of
$-a$, which implies that $x-x_0$ is the same multiple of $b$.

Let $w_{v}$ be the unique element of ${R_\omega}^2$ such that
$(w_{v})^\perp$ is the unique {\it solution} with $\frac{x_{w_{v}}}{a}
\in \D$.  As $x_{w_{v}}=-\,y_{(w_{v})^\perp}$, this is possible since,
by the above, the subset of $K_\omega$ consisting of the elements
$-\frac{y}{a}$, where $y$ varies over the second components of {\it
  solutions}, is exactly one orbit by translation under $R_\omega$
(without repetition).

Let us define $\ga_{v}=\begin{pmatrix} v & w_{v}\end{pmatrix}=
\begin{pmatrix} a & x_{w_{v}}\\ b & y_{w_{v}}\end{pmatrix}$. 
We have $\ga_v\in\Ga$ since ${w_{v}}^\perp$ is a {\it solution} so
that $\det \ga_{v}=1$. Furthermore $\ga_{v}\in \Ga^\sharp$ since $v\in
R_{\omega,\prim}^{2,\sharp}$.  Let $g= \ga_{v}$. By Equation
\eqref{eq:pppg}, the first column of $\ppp_g$ is
$(a\pi_\omega^{-\omega(a)},b\pi_\omega^{-\omega(a)})=
\pi_\omega^{\;\log_{q_\omega}|a|_\omega}\,v=\dirp{v}$, so that
$\ppp_g\in P^-_\Theta$ if and only if $\dirp{v}\in \Theta$. Since
$\norm{v}= |\,a\,|_\omega=q_\omega^{-\omega(a)}$ and by Equation
\eqref{eq:refineLU}, we have $\aaa_g\in A_{n}$ if and only if $\norm v
= q_\omega^{n}$. Again by Equation \eqref{eq:refineLU}, we have
$\uuu^+_g\in U^+_{\D'}$ if and only if $\frac{x_{w_{v}}}{x_{v}}
=\frac{x_{w_{v}}}{a} \in\D'$.

The map $v\mapsto \ga_{v}$ from $R_{\omega,\prim}^{2,\sharp}$ to
$\Ga^\sharp$ is clearly injective. Its image is $\Ga^\sharp\cap
(P^-\,U^+_\D)$, since if $\begin{pmatrix} v & w\end{pmatrix}
\in\Ga^\sharp\cap(P^-\,U^+_\D)$ and $v=(a,b)$, then $v$ belongs to
$R_{\omega,\prim}^{2,\sharp}$ and $w^\perp$ is a {\it solution} such
that by Equation \eqref{eq:refineLU} we have $-\frac{y_{w^\perp}}{a}
=\frac{x_{w}}{a} \in \D$, hence $w= w_{v}$ by uniqueness. We clearly
have $y_v=b\in I$ if and only if $\ga_v\in\Ga_0[I]$. This proves the
result.
\cqfd

\section{Joint equidistribution of primitive lattice points}
\label{sect:equidistrib}

The aim of this section is to prove the main result of this paper,
Theorem \ref{theo:main}, establishing the effective joint
equidistribution of directions and renormalized solutions to the
associated gcd equation for primitive lattice points, generalizing
Theorem \ref{theo:mainintro} in the Introduction to any function
field.

The main tool for this result is an adaptation of two theorems of
Gorodnik-Nevo \cite{GorNev12}, that we now state, after the necessary
definitions.

Let ${\bf G'}$ be an absolutely connected and simply connected
semi-simple algebraic group over $K_\omega$, which is almost
$K_\omega$-simple.  Let $G'={\bf G'}(K_\omega)$ be the locally compact
group of $K_\omega$-points of ${\bf G'}$. Let $\Ga'$ be a
non-uniform\footnote{This implies that ${\bf G'}$ is isotropic over
  $K_\omega$, as part of the assumptions of \cite{GorNev12}.}  lattice
in $G'$, and let $\mu_{G'}$ be any (left) Haar measure of $G'$. Note
that $G'=G$ and $\Ga'=\Ga_0[I]$ satisfy these assumptions for every
nonzero ideal $I$ of $R_\omega$.

Let $\rho>0$. Let $(\V'_\epsilon)_{\epsilon>0}$ be a fundamental
system of neighborhoods of the identity in $G'$, which 

$\bullet$~ is symmetric (that is, $x\in \V'_\epsilon$ if and only if
$x^{-1}\in \V'_\epsilon$),

$\bullet$~ is nonincreasing with $\epsilon$, and 

$\bullet$~ has {\it upper local dimension} $\rho$, that is, there
exist $m_1,\epsilon_1>0$ such that $\mu_{G'}(\V'_\epsilon)\geq
m_1\,\epsilon^\rho$ for every $\epsilon\in\;]0,\epsilon_1[\,$.

Let $C\geq 0$. Let $(\B_n)_{n\in\NN}$ be a family of measurable
subsets of $G'$.  We define
$$
(\B_n)^{+\epsilon}= \V'_\epsilon \B_n\V'_\epsilon=
\bigcup_{g,h\in\V'_\epsilon}g \B_n h \;\;\;{\rm and}\;\;\;
(\B_n)^{-\epsilon}= \bigcap_{g,h\in\V'_\epsilon}g \B_n h\;.
$$ 
The family $(\B_n)_{n\in\NN}$ is {\em $C$-Lipschitz well-rounded} with
respect to $(\V'_\epsilon)_{\epsilon>0}$ if there exists $\epsilon_0
>0$ and $n_0\in\NN$ such that for all $\epsilon\in\; ]0,\epsilon_0[$
and $n\geq n_0$, we have
$$
\mu_{G'}((\B_n)^{+\epsilon})\leq
(1+C\,\epsilon)\;\mu_{G'}((\B_n)^{-\epsilon})\;.
$$

\btheo\label{theo:GorodnikNevo} 
For every $\rho>0$, there exists $\tau(\Ga')\in \;]0,
\frac{1}{2(1+\rho)}]$ such that for every $C\geq 0$, for every
symmetric nonincreasing fundamental system
$(\V'_\epsilon)_{\epsilon>0}$ of neighborhoods of the identity in $G'$
with upper local dimension $\rho$, for every family $(\B_n)_{n\in\NN}$
of measurable subsets of $G'$ that is $C$-Lipschitz well-rounded with
respect to $(\V'_\epsilon)_{\epsilon>0}$, and for every $\delta>0$, we
have that, as $n\ra +\infty$,
$$
\Big|\card(\B_{n}\cap\Ga')-
\frac{1}{\|\mu_{\Ga'\bs G'}\|}\,\mu_{G'}(\B_{n})\Big|=
\bigO\big(\mu_{G'}(\B_{n})^{1-\tau(\Ga')+\delta}\big)\;,
$$
where the function $\bigO(\cdot)$ depends only on $G',\Ga',\delta,C,
(\V'_\epsilon)_{\epsilon>0},\rho$.
\etheo

\dem The proof is a simple adaptation of a particular case of results
of Gorodnik-Nevo \cite{GorNev12}, which are phrased for algebraic
number fields and not for function fields.

By the assumptions on ${\bf G'}$ and $\Ga'$, and by
\cite[Theo.~2.8]{AthGhoPra12}, the regular representation $\pi^0$ of
$G'$ on $\LL^2_0(G'/\Ga')$ has a spectral gap. By \cite{CowHaaHow88}
(see \cite[Theo.~2.7]{AthGhoPra12}), since $\pi^0$ has a spectral gap,
there exists $p\geq 2$ such that $\pi^0$ is strongly $\LL^p$ (called
$\LL^{p+}$ in \cite[Def.~3.1]{GorNev12}). We do not know what is the
smallest such $p$. As in \cite[Eq.~(3.1)]{GorNev12}, let $n_e(p)=1$ if
$p=2$ and $n_e(p)=\lceil\frac{p}{2}\rceil\in\NN-\{0,1\}$
otherwise. Since $\pi^0$ is strongly $\LL^p$, by
\cite[Theo.~4.5]{GorNev12}, for every measurable subset $B$ of $G'$
with finite and positive Haar measure, if $\pi^0(\beta)$ is the
operator on $\LL^2_0(G'/\Ga')$ defined by
$$
\pi^0(\beta)f(x)=\frac{1}{\mu_{G'}(B)}\, \int_{B}f(g^{-1}x)\,d\mu_{G'}(g)
$$ 
for all $f\in\LL^2_0(G'/\Ga')$ and almost all $x\in G'/\Ga'$, then we
have that, for every $\eta>0$,
$$
\|\pi^0(\beta)\|=
\bigO_{G',\Ga',\eta}\big((\mu_{G'}(B))^{-\frac{1}{2\,n_e(p)}+\eta}\big)\;.
$$

Actually, Theorem 4.5 of \cite{GorNev12} is stated in characteristic
zero. But its proof has two ingredients, a spectral transfer
principle, which is valid for any locally compact second countable
group by \cite[Theo.~1]{CowHaaHow88}, and a Kunze-Stein phenomenon,
which is valid even in positive characteristic by
\cite[Theo.~1]{Veca02}.

Now, by \cite[Theo.~1.9]{GorNev12} where $a=1$, which is valid for
any locally compact second countable group, and whose assumptions we
just verified, we have
$$
\Big|\frac{\card(\B_{n}\cap\Ga')}{\mu_{G'}(\B_{n})}-
\frac{1}{\|\mu_{\Ga'\bs G'}\|}\,\Big|=
\bigO_{G',\Ga',C,\rho,(\V'_\epsilon)_{\epsilon>0}} 
\big(\mu_{G'}(\B_{n})^{(-\frac{1}{2\,n_e(p)}+\eta)(\frac{1}{\rho+1})}\big)\;.
$$ 
Theorem \ref{theo:GorodnikNevo} follows with $\tau(\Ga')=
\frac{1}{2\,n_e(p)(\rho+1)}$.  
\cqfd

\bigskip
The main result that will allow us to use Theorem
\ref{theo:GorodnikNevo} is the following proposition. We will use, as
a fundamental system of neighborhoods of the identity element in $G$,
a family of compact-open subgroups of $G(\OOO_\omega)$ given by
the kernels of the morphisms of reduction modulo $\pi_\omega^n
\OOO_\omega$ for $n\in\NN$. For every $\epsilon>0$, let $N_\epsilon=
\big\lfloor-\log_{q_\omega} \epsilon\, \big\rfloor$ so that
$N_\epsilon\geq 1$ if and only if $\epsilon\leq \frac{1}{q_\omega}$. Let
$\V_\epsilon=G(\OOO_\omega)$ if $\epsilon>\frac{1}{q_\omega}$ and
otherwise let
\begin{align*}
\V_\epsilon&=
\ker(G(\OOO_\omega)\ra G(\OOO_\omega/\pi_\omega^{N_\epsilon}\OOO_\omega))
\\ &= \Big\{\begin{pmatrix}
  1+\pi_\omega^{N_\epsilon}\alpha &\pi_\omega^{N_\epsilon}\ga\\
  \pi_\omega^{N_\epsilon}\beta & 1+\pi_\omega^{N_\epsilon}\delta
  \end{pmatrix}
  \in G(\OOO_\omega):\alpha,\beta,\gamma,\delta \in\OOO_\omega\Big\} \;.
\end{align*}
The family $(\V_\epsilon)_{\epsilon>0}$ is indeed nondecreasing. Note
that for all $\epsilon_1,\dots, \epsilon_k>0$, we have
\begin{align*}
\min\{N_{\epsilon_1},\cdots, N_{\epsilon_k}\}&\geq
\min\{-\log_{q_\omega}\epsilon_1,\cdots, -\log_{q_\omega}\epsilon_k\}-1
\\ &\geq \min\{-\log_{q_\omega}(\epsilon_1+\cdots+\epsilon_k)\}-1
\geq N_{q_\omega(\epsilon_1+\cdots+\epsilon_k)}\,,
\end{align*}
hence
\begin{equation}\label{eq:minNepsilons}
  \V_{\epsilon_1}\V_{\epsilon_2}\cdots\V_{\epsilon_k}\;\;\subset\;\;
  \V_{q_\omega(\epsilon_1+\cdots+\epsilon_k)}\,.
\end{equation}

\bprop\label{prop:verifLWR} For all metric balls $\Theta$ in $\Snu$
and $\D'$ in $K_\omega$ with radius less than $1$, the family
$\big(P^-_{\Theta}\, A_{n}\, U^+_{\D'} \big)_{n\in\NN}$ is
$0$-Lipschitz well-rounded with respect to
$(\V_\epsilon)_{\epsilon>0}$.  
\eprop

\dem We will actually prove (as allowed by the ultrametric situation)
the stronger statement that given $\Theta$ and $\D'$ as above, if
$\epsilon$ is small enough, then for every $n\in\NN$,
$$
\big(P^-_\Theta A_n\,U^+_{\D'}\big)^{-\epsilon}=P^-_\Theta A_n\,U^+_{\D'}
=\big(P^-_\Theta A_n\,U^+_{\D'}\big)^{+\epsilon}\;.
$$

We start the proof by some elementary linear algebra considerations.
For every subgroup $H$ of $G$, let $\V_\epsilon^H = \V_\epsilon \cap
H$. We endow $\M_2(K_\omega)$ with its supremum norm
$\|\cdot\|_\omega$ defined, for every $X\in \M_2(K_\omega)-\{0\}$, by
$\|X\|_\omega= \max\{|X_{i,j}|_\omega:1\leq i,j\leq 2\} \in
q_\omega^\ZZ$. The unit ball of $\|\cdot\|_\omega$ is
$\M_2(\OOO_\omega)$. We denote the operator norm of a linear operator
$\ell$ of $\M_2(K_\omega)$ by
$$
\|\ell\|_\omega=\max\Big\{\frac{\|\ell(X)\|_\omega}{\|X\|_\omega} :
X\in\M_2(K_\omega)-\{0\}\Big\} \in q_\omega^\ZZ \cup\{0\}\,,
$$
so that $\ell(\M_2(\OOO_\omega))\subset
\M_2(\pi_\omega^{-\log_{q_\omega}\|\ell\|_\omega}\OOO_\omega)$.  For
every $g\in G$, recall that $\operatorname{Ad} g$ is the linear
automorphism $x\mapsto gxg^{-1}$ of $\M_2(K_\omega)$.

\blemm\label{eq:calVdivgroup}
For all $\epsilon>0$ and $g\in G$, we have
$$
g\,\V_\epsilon\, g^{-1}\;\subset\;
\V_{\epsilon\,\|\operatorname{Ad} g\,\|_\omega}\;,\;\;\;
\V_\epsilon=\V^{P^-}_\epsilon\;\V^{U^+}_\epsilon\;\;\;{\rm and}\;\;\;
\V^{P^-}_\epsilon=\V^{U^-}_\epsilon\;\V^{Z}_\epsilon\;.
$$
Furthermore, we have $\mu_G(\V_\epsilon)\geq
\frac{q_\omega^{\,2}}{q_\omega^{\,2}-1}\, \epsilon^3$ for every
$\epsilon>0$ small enough, so that $\rho=3$ is an upper local
dimension of the family $(\V_\epsilon)_{\epsilon>0}$.
\elemm

\dem Let $I_2$ be the identity element in $G$. The first claim follows from
the fact that
$$
g\,\V_\epsilon \,g^{-1}=
I_2+\pi_\omega^{N_\epsilon}g\M_2(\OOO_\omega)g^{-1}
\subset
I_2+\pi_\omega^{N_\epsilon-\log_{q_\omega}\|\operatorname{Ad} g\,\|_\omega} \M_2(\OOO_\omega)=
\V_{\epsilon\,\|\operatorname{Ad} g\,\|_\omega}\;.
$$
The second and third claims follow from the fact that by Equations
\eqref{eq:refineLU} and \eqref{eq:pppg}, if $g\in \V_\epsilon$ then
$\aaa_g=I_2$, $\uuu^\pm_g\in \V^{U^\pm}_\epsilon$ and $\mmm_g\in
\V^{Z}_\epsilon$.

Let us now apply Lemma \ref{lem:decompHaarG} with the decomposition
$\V_\epsilon=\V^{U^-}_\epsilon\;\V^{Z}_\epsilon\; \V^{U^+}_\epsilon$:
$$
\mu_G(\V_\epsilon)=\frac{q_\omega}{q_\omega+1}\,
\mu_{U^-}(\V^{U^-}_\epsilon)\;\mu_{Z(\OOO_\omega)}(\V^{Z}_\epsilon)\;
\mu_{U^+}(\V^{U^+}_\epsilon)\;.
$$
By Equation \eqref{eq:muUpm} applied twice, by the left part of
Equation \eqref{eq:haarKmodR}, and since $N_\epsilon=\lfloor
-\log_{q_\omega} \epsilon\,\rfloor$, we have that for $\epsilon\leq
\frac{1}{q_\omega}$,
\begin{align*}
\mu_G(\V_\epsilon)& =\frac{q_\omega}{q_\omega+1}\,
\mu_{K_\omega}(\pi_\omega^{N_\epsilon}\OOO_\omega)\;
\big|\OOO_\omega^\times/(1+\pi_\omega^{N_\epsilon}\OOO_\omega)\big|^{-1}\;
\mu_{K_\omega}(\pi_\omega^{N_\epsilon}\OOO_\omega)
\\ & =\frac{q_\omega}{q_\omega+1}\;
\frac{1}{(q_\omega-1)\,q_\omega^{N_\epsilon-1}}\;q_\omega^{-2\,N_\epsilon}
\geq \frac{q_\omega^{\,2}}{q_\omega^{\,2}-1}\,\epsilon^3\;.
\end{align*}
This proves the final claim of Lemma \ref{eq:calVdivgroup}.
\cqfd

\medskip
The main ingredient in the proof of Proposition \ref{prop:verifLWR} is
the following effective refined LU decomposition.

\blemm\label{eq:effectiveLU} With $c:G\ra\;]0,+\infty[$ the
continuous function defined by $h\mapsto\|\operatorname{Ad}
h\,\|_\omega$, for every $g\in G$ with $|\chi(\aaa_g)|_\omega\leq 1$, we
have
$$
\V_\epsilon \;g\;\V_{\epsilon}
\;\subset\; \ppp^-_g\;\V^{P^-}_{q_\omega(c(\ppp_g)+c(\uuu_g))\epsilon}\;
\aaa_g\;\V^{U^+}_{q_\omega(c(\ppp_g)+2c(\uuu_g))\epsilon}\;\uuu_g\;.
$$ 
\elemm

\dem In order to simplify notation, let $\aaa=\aaa_g$, $\ppp=\ppp^-_g$
and $\uuu=\uuu^+_g$, so that $g=\ppp\,\aaa\,\uuu$. For every $h\in G$,
let $c_h= \|\operatorname{Ad} h\,\|_\omega$.  In the following sequence
of equalities and inclusions, we use

$\bullet$~ 
the first claim of Lemma \ref{eq:calVdivgroup}, for the first
inclusion,

$\bullet$~ the second claim of Lemma \ref{eq:calVdivgroup}, for the
second equality,

$\bullet$~ the fact that 
$$
\aaa\, \V^{P^-}_{c_\uuu\epsilon}= \aaa\,
\V^{U^-}_{c_\uuu\epsilon}\;\V^{Z}_{c_\uuu\epsilon}\subset
\V^{U^-}_{c_\uuu\epsilon}\,\aaa\,\V^{Z}_{c_\uuu\epsilon}=\V^{P^-}_{c_\uuu\epsilon}
\,\aaa\subset \V_{c_\uuu\epsilon}\,\aaa
$$ 
by the third claim of Lemma \ref{eq:calVdivgroup}, by the
right-hand side of Equation \eqref{eq:dilatlarate} with
$\chi(\aaa)\in\OOO_\omega$ and since $\aaa$ and $Z$ commute, for the
second inclusion,

$\bullet$~ the facts that $\V_{c_\uuu\epsilon}$ is a normal subgroup
of $G(\OOO_\omega)$ and that $\V^{U^+}_{c_\ppp\epsilon}\subset
G(\OOO_\omega)$, for the third equality,

$\bullet$~ again the second claim of Lemma \ref{eq:calVdivgroup}, and
the left-hand side of Equation \eqref{eq:dilatlarate} with
$\chi(\aaa)\in\OOO_\omega$, for the last inclusion.

\noindent We thus have
\begin{align*}
\V_\epsilon \,g\,\V_\epsilon& =
\ppp\,\ppp^{-1}\V_\epsilon \,\ppp\,\aaa\,\uuu\,\V_\epsilon \,\uuu^{-1}\uuu
\;\subset\; \ppp\,\V_{c_\ppp\epsilon}\,\aaa\,\V_{c_\uuu\epsilon} \,\uuu
\;=\; \ppp\,\V^{P^-}_{c_\ppp\epsilon}\,\V^{U^+}_{c_\ppp\epsilon}\,\aaa\,
\V^{P^-}_{c_\uuu\epsilon}\V^{U^+}_{c_\uuu\epsilon} \,\uuu \\ &\subset\;
\ppp\,\V^{P^-}_{c_\ppp\epsilon}\,\V^{U^+}_{c_\ppp\epsilon}\,
\V_{c_\uuu\epsilon}\,\aaa\,\V^{U^+}_{c_\uuu\epsilon}  \,\uuu \;=\;
\ppp\,\V^{P^-}_{c_\ppp\epsilon}\,\V_{c_\uuu\epsilon}\,
\V^{U^+}_{c_\ppp\epsilon}\,\aaa\,\V^{U^+}_{c_\uuu\epsilon}  \,\uuu\\ &
\subset\;\ppp\,\V^{P^-}_{c_\ppp\epsilon}\,\V^{P^-}_{c_\uuu\epsilon}\,\aaa\,
\V^{U^+}_{c_\uuu\epsilon}\V^{U^+}_{c_\ppp\epsilon}\,\V^{U^+}_{c_\uuu\epsilon} \,\uuu
\;.
\end{align*}
Lemma \ref{eq:effectiveLU} now follows from Equation \eqref{eq:minNepsilons}.  
\cqfd

\medskip
Now, in order to prove Proposition \ref{prop:verifLWR}, we write
$\Theta=v_0+\pi_\omega^{m} {\OOO_\omega}^2$ and $\D'=x_0+\pi_\omega^{m'}
\OOO_\omega$, for some $m,m'\in\NN-\{0\}$, $x_0\in K_\omega$ and $v_0\in
\Snu$.  Let
$$
c=\max\{q_\omega(c(\ppp)+2c(\uuu)) : \ppp\in P^-_\Theta,\;\uuu\in
U^+_{\D'}\}\,,
$$
which is finite since $P^-_\Theta$ and $U^+_{\D'}$ are compact. Let
$\epsilon_0=\frac{1}{c}\,q_\omega^{-m'-m}>0$, so that we have
$N_{c\epsilon}> \max\{m,m'\}\geq 1$ if $\epsilon<\epsilon_0$.

Let us fix $\epsilon\in\;]0,\epsilon_0[$. We claim that
\begin{equation}\label{eq:petitpertubpaspertub}
P^-_\Theta \,\V^{P^-}_{c\epsilon}=P^-_\Theta \;\;\;{\rm and}\;\;\;
\V^{U^+}_{c\epsilon}\,U^+_{\D'}=U^+_{\D'}\;.
\end{equation}
Indeed, the inclusion of the right-hand sides into the left-hand sides
of these equalities are immediate. If $\ppp\in P^-_\Theta$ and
$\ppp'\in \V^{P^-}_{c\epsilon}$, we may write
$$
\ppp=\begin{pmatrix} x_{v_0}+\pi_\omega^{m}\alpha & 0 \\
y_{v_0}+\pi_\omega^{m}\beta & (x_{v_0}+\pi_\omega^{m}\alpha)^{-1}\end{pmatrix} 
\;\;\;{\rm and}\;\;\;
\ppp'=\begin{pmatrix} 1+\pi_\omega^{N_{c\epsilon}}\alpha' & 0 \\
\pi_\omega^{N_{c\epsilon}}\beta' & (1+\pi_\omega^{N_{c\epsilon}}\alpha')^{-1}
\end{pmatrix} 
$$
for some $\alpha,\beta,\alpha',\beta'\in\OOO_\omega$, so that
$$
\ppp\,\ppp'=\begin{pmatrix} 
x_{v_0}+\pi_\omega^{m}\alpha + \pi_\omega^{N_{c\epsilon}}\alpha''& 0 \\
y_{v_0}+\pi_\omega^{m}\beta + \pi_\omega^{N_{c\epsilon}}\beta''& 
(x_{v_0}+\pi_\omega^{m}\alpha+\pi_\omega^{N_{c\epsilon}}\alpha'')^{-1}\end{pmatrix}
$$ 
for some $\alpha'',\beta''\in\OOO_\omega$ (since $x_{v_0},y_{v_0}
\in\OOO_\omega$). The first claim then follows from the fact that
$N_{c\epsilon}>m$.  The inclusion $\V^{U^+}_{c\epsilon}\,U^+_{\D'}
\subset U^+_{\D'}$ follows from a similar and even easier computation.

Now for every $n\in\NN$, we have by Lemma \ref{eq:effectiveLU} and
Equation \eqref{eq:petitpertubpaspertub} that
$$ 
\big(P^-_\Theta A_n\,U^+_{\D'}\big)^{+\epsilon}=\V_\epsilon\,
P^-_\Theta A_n\,U^+_{\D'}\V_\epsilon\subset P^-_\Theta
\,\V^{P^-}_{c\epsilon} A_n\,\V^{U^+}_{c\epsilon}\,U^+_{\D'}=
P^-_\Theta A_n\,U^+_{\D'}\;.
$$ 
Since the converse inclusion is immediate, we have $\big(P^-_\Theta
A_n\,U^+_{\D'}\big)^{+\epsilon}=P^-_\Theta A_n\,U^+_{\D'}$. This
implies that $g\, P^-_\Theta A_n\,U^+_{\D'}\;h = P^-_\Theta
A_n\,U^+_{\D'}$ for all $g,h\in \V_\epsilon$ so that $ \big(P^-_\Theta
A_n\,U^+_{\D'}\big)^{-\epsilon} = P^-_\Theta A_n\,U^+_{\D'}$.  

This concludes the proof of Proposition \ref{prop:verifLWR}.
\cqfd

\bigskip
The main result of this paper is the following one. Recall that $z_v$,
$z'_v$ and $\dirp v$ for $v$ in ${K_\omega}^2-\{0\}$ have been defined
in Equations \eqref{eq:znu} and \eqref{eq:dirpv}. If $v=(a,b)\in
R^2_{\omega,\prim}$, we denote by $w_v$ any element of
$R^2_{\omega,\prim}$ such that $(w_v)^\perp=(x,y)$ is a solution to
the equation $ax+by=1$. As seen in the proof of Proposition
\ref{prop:matrepprimvect} if $|\,a\,|_\omega\geq |\,b\,|_\omega$,
and by symmetry otherwise, the class $\frac{z_{w_{v}}}{z_v}+R_\omega$
of $\frac{z_{w_{v}}}{z_v}$ in the quotient $K_\omega/R_\omega$ does
not depend on the choice of $w_v$. For every nonzero ideal $I$ of
$R_\omega$, let
\begin{equation}\label{eq:defcsubI}
 c_I= \frac{(q_\omega^{\,2}-1)\,\zeta_K(-1)\,N(I)\,\prod_{\ppp | I}
  \big(1+\frac{1}{N(\ppp)}\big)}{q_\omega^{2}}
\end{equation}

\btheo\label{theo:main} For every nonzero ideal $I$ of $R_\omega$, for
the weak-star convergence on the compact space $\Snu\times
(K_\omega/R_\omega)$, we have, as $n\ra+\infty$,
$$
 c_I\;q_\omega^{-2n}
\sum_{v\in R^2_{\omega,\prim}\;:\;\|v\|_\omega=q_\omega^n,\;z'_v\in I} 
\Delta_{\dirp v}\otimes\Delta_{\frac{z_{w_{v}}}{z_v}+R_\omega}\;\;\weakstar\;\;
\mu_{\Snu}\otimes \mu_{K_\omega/R_\omega}\;.
$$ 
Furthermore, there exists $\tau\in\;]0,\frac{1}{8}]$ such that for all
$\epsilon,\delta>0$, there is a multiplicative error term in the above
equidistribution claim of the form
$1+\bigO_{\omega,\delta,I}\big(q_\omega^{2n(-\tau+\delta)}
\,\|f\|_\epsilon\,\|g\|_\epsilon\big)$ when evaluated on $(f,g)$ for
all $\epsilon$-locally constant maps $f:\Snu \ra\RR$ and
$g:K_\omega/R_\omega\ra\RR$:
\begin{align*}
  c_I&\;q_\omega^{-2n} 
\sum_{v\in R^2_{\omega,\prim}\;:\;\|v\|_\omega=q_\omega^n,\;z'_v\in I} 
f(\,\dirp{v}\,)\;g\big(\frac{z_{w_{v}}}{z_v}+R_\omega\big)\\ & =
\Big(\int_{\Snu} f\,d\mu_{\Snu}\Big)
\Big(\int_{K_\omega/R_\omega}g\,d\mu_{K_\omega/R_\omega}\Big)\Big(1+\bigO_{\omega,\delta}
\big(q_\omega^{2n(-\tau+\delta)} \,\|f\|_\epsilon\,\|g\|_\epsilon\big)\Big)\;.
\end{align*}
\etheo

When ${\bf C}=\PP^1$, $\omega= \omega_\infty$ and $I=
R_{\omega_\infty}$, we recover Theorem \ref{theo:mainintro} in the
Introduction by using Equations \eqref{eq:defcsubI},
\eqref{eq:indexHecke} and \eqref{zetamoinsun}, as well as the fact
that $q_\omega=q$. Note that up to changing the constant $C_I$, the
same result holds when $v$ ranges over the elements in
$R^2_{\omega,\prim}$ with $|v\|_\omega\leq q_\omega^n$ and $z'_v\in
I$. Also note that in the statement of Theorem \ref{theo:main}, the
measures $\mu_{\Snu}$ and $\mu_{K_\omega/R_\omega}$ are not normalized
to be probability measures, see Equations \eqref{eq:mesSnu} and
\eqref{eq:haarKmodR} if a normalization is useful, as for instance in
Corollary \ref{coro:counting}.

Given a nonzero (possibly nonprincipal) ideal $J$ of $R_\omega$, an
effective joint equidistribution result similar to the one of Theorem
\ref{theo:main} is possible when the elements $v=(a,b)\in
{R_{\omega}}^2$ are not assumed to be primitive, but to satisfy
that $a,b$ generate the ideal $J$.

\medskip
\dem 
Let $I$ be a nonzero ideal of $R_\omega$. Let $\tau= \tau(\Ga_0[I])
\in\;]0,\frac{1}{8}[$ be as in Theorem \ref{theo:GorodnikNevo} applied
with $G'=G$ and $\Ga'=\Ga_0[I]$, and with $(\V'_\epsilon)_{\epsilon>0}
= (\V_\epsilon)_{\epsilon>0}$ which has upper local dimension $\rho=3$
according to the final claim of Lemma \ref{eq:calVdivgroup}.
Let $\delta\in\;]0,\tau]$. Fix a compact-open strict 
fundamental domain $\D$ for the action by translations of $R_\omega$
on $K_\omega$, such that for all $x_0\in \D$ and $m'\in\NN-\{0\}$, we
have $B(x_0,q_\omega^{-{m'}})= x_0+\pi_\omega^{m'} \OOO_\omega \subset
\D$. This is possible since $R_\omega \cap \pi_\omega \OOO_\omega =
\{0\}$ by Equation \eqref{eq:interRnupiOnu}.  Note that for all
$v_0\in \Snu$ (respectively $v_0\in \Snusharp$) and $m\in \NN-\{0\}$,
the ball $B(v_0,q_\omega^{-m})=v_0+\pi_\omega^{m} {\OOO_\omega}^2$ is
contained in $\Snu$ (respectively $\Snusharp$).

Let us prove that for all $m,m'\in\NN-\{0\}$,  $x_0\in \D$ and
$v_0\in \Snu$, if $\Theta=v_0+\pi_\omega^{m} {\OOO_\omega}^2$ and $\D'=x_0+
\pi_\omega^{m'}\OOO_\omega$, then, as $n\ra+\infty$
\begin{align}
\card&\big\{v\in R^2_{\omega,\prim}\;:\; \|v\|_\omega=q_\omega^n,\;z'_v\in I,
\;\dirp{v}\in\Theta,\;\frac{z_{w_v}}{z_v}\in \D'\big\}\nonumber
\\ &=\label{eq:reduc1} \frac{1}{c_I}
\;q_\omega^{2n}\,\mu_{\Snu}(\Theta)\; \mu_{K_\omega}(\D')\big(1+
\bigO_{\omega,\delta,I}\big(q_\omega^{2n(-\tau+\delta)} q_\omega^{m+m'}\big)\big)\;.
\end{align}
Since the characteristic functions $\mathbbm{1}_\Theta$ and
$\mathbbm{1}_{\D'}$ of $\Theta$ and $\D'$ are respectively
$q_\omega^{-m}$- and $q_\omega^{-m'}$-locally constant, and by a
finite additivity argument, this proves Theorem \ref{theo:main}.

We first claim that in order to prove the counting result of elements
in $R^2_{\omega,\prim}$ stated in Equation \eqref{eq:reduc1}, we only
have to prove an analogous counting result of elements in
$R^{2,\sharp}_{\omega,\prim}$, namely that for all $m,m'\in\NN-\{0\}$,
for all $x_0\in \D$ and $v_0\in \Snusharp$, if
$\Theta=v_0+\pi_\omega^{m} {\OOO_\omega}^2$ and $\D'=x_0+
\pi_\omega^{m'}\OOO_\omega$, then, as $n\ra+\infty$
\begin{align}
\card&\big\{v\in R^{2,\sharp}_{\omega,\prim}\;:\; \|v\|_\omega= q_\omega^n,
\;y_v\in I,\;\dirp{v}\in\Theta,\;\frac{x_{w_v}}{x_v}\in \D'\big\}
\nonumber \\ &=\label{eq:reduc2}
\frac{1}{c_I}
\;q_\omega^{2n}\,\mu_{\Snu}(\Theta)\; \mu_{K_\omega}(\D')\big(1
+\bigO_{\omega,\delta,I}\big(q_\omega^{2n(-\tau+\delta)} q_\omega^{m+m'}\big)\big)\;.
\end{align}
Indeed, by Lemma \ref{lem:maxmin} and Equation \eqref{eq:znu}, since
$\det\begin{pmatrix}v& w_v\end{pmatrix}=1$, we have
$\frac{z_{w_v}}{z_v}=\frac{x_{w_v}}{x_v}$ when $v$ belongs to
$R^{2,\sharp}_{\omega,\prim}$ except finitely many of them.  The
involutive linear map $\iota= \big(\begin{smallmatrix}0 & 1\\ 1 &
  0\end{smallmatrix}\big)$ of exchange of coordinates

$\bullet$~ preserves the subsets $R^{2}_{\omega,\prim}$ and $\Snu$ of
  the plane ${K_{\omega}}^2$,

$\bullet$~ sends the compact-open set $\Snu-\Snusharp$ into
  $\Snusharp$,

$\bullet$~ sends an element $v$ in $R^{2}_{\omega,\prim}-
  R^{2,\sharp}_{\omega,\prim}$ to the element $\iota(v)$ in
  $R^{2,\sharp}_{\omega,\prim}$ such that $z'_v=z_{\iota(v)}=
  y_{\iota(v)}$ and $\frac{z_{w_v}}{z_v}
  =\frac{x_{w_{\iota(v)}}}{x_{\iota(v)}}$ again by Lemma
  \ref{lem:maxmin} and Equation \eqref{eq:znu}, and

$\bullet$~ sends $v_0+\pi_\omega^{m} {\OOO_\omega}^2$ to
$\iota(v_0)+\pi_\omega^{m} {\OOO_\omega}^2$.

\noindent Hence Equation \eqref{eq:reduc1} follows from Equation
\eqref{eq:reduc2}

\medskip
Now according to Proposition \ref{prop:verifLWR}, the family
$\big(P^-_{\Theta}\, A_{n}\, U^+_{\D'} \big)_{n\in\NN}$ is
$0$-Lipschitz well-rounded in $G$ with respect to
$(\V_\epsilon)_{\epsilon>0}$. Note that
$$
\Ga\cap (P^-_{\Theta}\, A_{n}\, U^+_{\D'})= 
\Ga^\sharp\cap (P^-_{\Theta}\, A_{n}\, U^+_{\D'})
$$
since $\Theta$ is contained in $\Snusharp$.  In the following
sequence of equalities, we use respectively

\smallskip
$\bullet$~ Proposition \ref{prop:matrepprimvect},

$\bullet$~ Theorem \ref{theo:GorodnikNevo} applied with $G'=G$,
$\Ga'=\Ga_0[I]$ and $(\B_n)_{n\in\NN}=(P^-_{\Theta}\, A_{n}\,
U^+_{\D'})_{n\in\NN}$,

$\bullet$~ Equation \eqref{eq:calchaarPThetaA0nUDom},

$\bullet$~ the fact that $\Theta$ is a metric ball of radius
$q_\omega^{-m}$ in the plane ${K_{\omega}}^2$ and $\D'$ a metric ball
of radius $q_\omega^{m'}$ in the line $K_{\omega}$.

\smallskip
\noindent We thus have
\begin{align}
\card&\big\{v\in R^{2,\sharp}_{\omega,\prim}\;:\; \|v\|_\omega= q_\omega^n,
\;y_v\in I,\;\dirp{v}\in\Theta,\;\frac{x_{w_v}}{x_v}\in \D'\big\}
\nonumber \\ &=\nonumber
\card\big(\Ga_0[I]\cap(P^-_{\Theta}\, A_{n}\, U^+_{\D'})\big)\\ &= 
\frac{\mu_G(P^-_{\Theta}\, A_{n}\, U^+_{\D'})}{\|\mu_{\Ga_0[I]\bs G}\|}
+ \bigO_{\omega,\delta,I}\big(\big(\mu_G(P^-_{\Theta}\, A_{n}\, U^+_{\D'})
\big)^{1-\tau+\delta}\big) \nonumber\\ &=
\frac{q_\omega^{2n+2}}{(q_\omega^{\,2}-1)\,\|\mu_{\Ga_0[I]\bs G}\|}\;
\mu_{\Snu}(\Theta)\;\mu_{K_\omega}(\D')+\bigO_{\omega,\delta,I}\big(
\big(q_\omega^{2n}\mu_{\Snu}(\Theta)\;\mu_{K_\omega}(\D')\big)^{1-\tau+\delta}
\big) \nonumber\\ &=
\frac{q_\omega^{2n+2}}{(q_\omega^{\,2}-1)\,\|\mu_{\Ga_0[I]\bs G}\|}\;
\mu_{\Snu}(\Theta)\;\mu_{K_\omega}(\D')\big(1+\bigO_{\omega,\delta,I}\big(
q_\omega^{2n(-\tau+\delta)}q_\omega^{2m(\tau-\delta)}q_\omega^{m'(\tau-\delta)}\big)\big)
\;.\label{eq:deducGN}
\end{align}
Since by Equations \eqref{eq:covolGa} and
\eqref{eq:indexHecke}, we have
$$
\|\mu_{\Ga_0[I]\bs G}\|=\|\mu_{\Ga\bs G}\|\;[\,\Ga:\Ga_0[I]\,]=
\zeta_K(-1)\, N(I)\,\prod_{\ppp | I}\big(1+\frac{1}{N(\ppp)}\big)
$$ 
and since $\tau\leq \frac{1}{8}$, this proves Equation
\eqref{eq:reduc2} and completes the proof of Theorem \ref{theo:main}.
\cqfd

\medskip 
We conclude this section by stating a counting result, which follows
from the equidistribution claim of Theorem \ref{theo:main} by
integrating on the pairs of constant functions with value $1$ on
$\Snu$ and on $K_\omega/R_\omega$, and by using Equations
\eqref{eq:mesSnu} and \eqref{eq:haarKmodR}.

\bcoro\label{coro:counting}
There exists  $\tau\in\;]0,\frac{1}{8}]$ such that for every
$\delta>0$, we have
\begin{align*}
  \card\;&\{v\in R^2_{\omega,\prim}\;:\;\|v\|_\omega=q_\omega^n,\;z'_v\in I\}
  \\ &= \frac{q^{{\bf g}-1}}{\zeta_K(-1)\;N(I)\,
  \prod_{\ppp | I}\big(1+\frac{1}{N(\ppp)}\big)} \;q_\omega^{2n} +
\bigO_{\delta,I}\big(q_\omega^{2n(1-\tau+\delta)}\big)\;.\;\;\;\Box
\end{align*}
\ecoro

\section{Application to the distribution of continued fraction
  expansions}
\label{sect:cfe}

In this section, we assume that ${\bf C}=\PP^1$ and $\omega=
\omega_\infty$, so that the notation in Section
\ref{subsect:functionfields} coincides with the notation of the
introduction: $K=\FF_q(Y)$, $R_{\omega_\infty}= R= \FF_q[Y]$,
$K_{\omega_\infty}= \wh K= \FF_q((Y^{-1}))$, $\OOO_{\omega_\infty} =
\OOO= \FF_q[[Y^{-1}]]$ and $|\cdot|_{\omega_\infty}=|\cdot|$.

Let us recall elementary facts on the continued fraction expansions in
$\wh K$, similar to the ones in $\RR$, see for instance the surveys
\cite{Lasjaunias00,Schmidt00}, and \cite{Paulin02} for a geometric
interpretation.  Any element $f\in \wh{K}$ may be uniquely written
$f=[f]+\{f\}$ with $[f]\in R$ (called the {\it integral part} of $f$)
and $\{f\}\in Y^{-1}\OOO$ (called the {\it fractional part} of $f$).
The {\it Artin map} $\Psi: Y^{-1}\OOO-\{0\}\ra Y^{-1}\OOO$ is defined
by $f\mapsto \big\{\frac{1}{f}\big\}$.  Any $f\in K-R$ has a unique
finite continued fraction expansion
$$
f=a_0 +
\cfrac{1}{a_1+\cfrac{1}{a_2+ \cfrac{1}{\cdots +\cfrac{1}{a_n}}}}\;,
$$ 
with $a_0=[f]\in R$ and $a_i= \big[\frac{1}{\Psi^{i-1}(f-a_0)} \big]$
a nonconstant polynomial for $1\leq i\leq n $ (called the {\it
  coefficients} of the continued fraction expansion of $f$), where
$n\in\NN-\{0\}$ is such that $\Psi^{n}(f-a_0)=0$.

Two finite sequences of polynomials $(P_i)_{-1\leq i\leq n}$ and
$(Q_i)_{-1\leq i\leq n}$ in $R$ are defined inductively as follows
$$
\begin{array}{cccc}
P_{-1}=1 & P_{0}=a_{0}, &   & P_{i}=P_{i-1}a_{i}+P_{i-2}\\
Q_{-1}=0 & Q_{0}=1, &   & Q_{i}=Q_{i-1}a_{i}+Q_{i-2}
\end{array}
$$
for $1\leq i\leq n$. The elements $P_i/Q_i$ for $0\leq i\leq n-1$ are
called the {\it convergents} of $f$, and $P_n/Q_n=f$. The convergents
have the following characterisation (see for instance
\cite[p.~140]{Schmidt00}): for all $P,Q\in R$ such that $\deg Q<
\deg Q_n$
\begin{equation}\label{eq:caracconv}
  {\rm if~} |f-P/Q|<\frac{1}{|Q|^2}\;\;{\rm then}
  \;\;P/Q \;\;{\rm is~a~convergent.}
\end{equation}
For $0\leq i\leq n-1$, we have
\begin{equation}\label{eq:approxconv}
\Big|\;f-\frac{P_i}{Q_i}\;\Big|=\frac{1}{|Q_i|\;|Q_{i+1}|}
\end{equation}
by for instance \cite[Eq.~(1.12)]{Schmidt00}), and
\begin{equation}\label{eq:PQcoprime}
Q_{i+1}P_{i}-P_{i+1}Q_{i} =(-1)^{i+1}\;.
\end{equation}
Since $\deg a_i\geq 1$ if $i\geq 1$, we have $\deg Q_i >\deg Q_{i-1}$
for $1\leq i\leq n$. If $f\in Y^{-1}\OOO$, then $a_0=0$ and
$P_i/Q_i\in Y^{-1}\OOO$, or equivalently $|P_i|<|Q_i|$, for $1\leq
i\leq n$.

The following result relates the shortest solutions to an equation
$ax+by=1$ with the continued fraction expansion of $a/b$.

\blemm \label{lem:shortestvsscf}
Let $a,b\in R-\{0\}$ be two coprime polynomials such that $a/b \in
\pi^{-1}\OOO$. Let $(P_i/Q_i)_{0\leq i\leq n}$ be the sequence of
convergents of $a/b$.  Then there exists a unique
$\lambda\in\FF_q^\times$ such that $(a,b)=(\lambda \,P_{n},\lambda\,
Q_{n})$ and $(-(-1)^n\lambda\, Q_{n-1},(-1)^n\lambda\, P_{n-1})$ is
the unique shortest solution to the equation $ax+by=1$.
\elemm

Note that this result implies that for all $a,b\in R-\{0\}$, the
equation $ax+by=1$ has one and only one shortest solution, up to
exchanging $a$ and $b$ if $|a|> |b|$ and to replacing $(a,b)$ by
$(a-\lambda' b, b)$ for the unique $\lambda'\in \FF_q^\times$ such that
$\deg (a-\lambda' b)<\deg b$ if $|a|= |b|$.

\medskip
\dem We may assume that $a\notin \FF_q^\times$, otherwise the result
is immediate with $\lambda=-a^{-1}$ since $P_0=0$ and $Q_0=1$. 

Since $P_n$ and $Q_n$ are coprime polynomials by Equation
\eqref{eq:PQcoprime} and $P_n/Q_n= a/b$, there exists $\lambda\in
\FF_q^\times$ such that $a=\lambda P_n$ and $b=\lambda Q_n$. Let
$a'=\frac{a}{(-1)^n\lambda}$ and $b'=\frac{a}{(-1)^n\lambda}$.  By
{\it solution}, we now understand a solution to the equation
$a'x+b'y=1$ with unknown $(x,y)\in R^2$.

We have $a'=(-1)^n P_n$ and $b'=(-1)^n Q_n$. Again by Equation
\eqref{eq:PQcoprime}, this implies that $(-Q_{n-1},P_{n-1})$ is a {\it
  solution}.

Let $(x'_0,y'_0)$ be another {\it solution}. Since we have
$|a'|<|b'|$, it follows from Lemma \ref{lem:maxmin} that $|x'_0|\geq
|y'_0|$, so that $\|(x'_0,y'_0)\|_\infty=|x'_0|$. We have
$\|(-Q_{n-1},P_{n-1})\|_\infty= |Q_{n-1}|$ since $P_{n-1}/Q_{n-1}\in
Y^{-1}\OOO$. In order to prove that $(-Q_{n-1},P_{n-1})$ is the unique
shortest {\it solution}, let us assume that $|x'_0|\leq |Q_{n-1}|$,
and let us prove that $(x'_0,y'_0)=(-Q_{n-1},P_{n-1})$.

Since $|x'_0|\leq |Q_{n-1}|<|Q_{n}|$, we have
$$
\Big|\;\frac{y'_0}{-x'_0}-\frac{P_n}{Q_n}\;\Big|=\frac{1}{|x'_0|\;|Q_{n}|}
<\frac{1}{|-x'_0|^2}\;.
$$
Hence by Equation \eqref{eq:caracconv}, $\frac{y'_0}{-x'_0}$ is a
convergent of $\frac{P_n}{Q_n}$, that is, there exists $i\in\{0,\dots,
n-1\}$ such that $\frac{y'_0}{-x'_0}=\frac{P_i}{Q_i}$. This implies in
particular that there exists $\lambda'\in\FF_q^\times$ such that
$(y'_0,-x'_0)=(\lambda'P_i,\lambda'Q_i)$.  Using Equation
\eqref{eq:approxconv} for the last equality, we have
$$
\frac{1}{|Q_{i}|\;|Q_{n}|}= \frac{1}{|x'_0|\;|Q_{n}|}=
\Big|\;\frac{y'_0}{-x'_0}-\frac{P_n}{Q_n}\;\Big|=
\Big|\;\frac{P_i}{Q_i}-\frac{P_n}{Q_n}\;\Big|
=\frac{1}{|Q_i|\;|Q_{i+1}|}\;.
$$
Since $|Q_{i+1}|<|Q_{n}|$ if $i<n-1$, this implies that $i=n-1$. Since
$(x'_0,y'_0)$ is a {\it solution}, we have $\lambda'=1$. Hence
$(y'_0,-x'_0) =(P_{n-1},Q_{n-1})$ as wanted.

Since the pair $(x_0,y_0)$ is a solution to the equation $ax+by=1$
if and only if the pair $((-1)^n\lambda\, x_0, (-1)^n\lambda \,y_0)$ is a
solution to the equation $a'x+b'y=1$, the result follows.
\cqfd

\medskip
The following result is an analogue in the field of formal Laurent
series to the main result of \cite{DinSin90} in the real field. It
gives an application of Theorem \ref{theo:mainintro} to the
distribution properties of the continued fraction expansions of
elements of $K$. For every $v=(a,b)\in R^2_\prim$, we denote by
$\Big(\frac{P_i(v)}{Q_i(v)}\Big)_{1\leq i\leq n_v}$ the continued
fraction expansion of $\frac{a}{b}$. We denote by $\mu_{Y^{-1}\OOO}$
the Haar measure of the compact additive group $Y^{-1}\OOO$,
normalized to be a probability measure.

\bcoro\label{coro:cfe} Let $P'=\prod_{i=1}^k\pi_i$ be a nonzero
polynomial in $R$, with prime factors $\pi_1,\dots,$ $\pi_k$. For the
weak-star convergence of measures on $Y^{-1}\OOO$, we have, as
$n\ra+\infty$,
$$
\frac{q^{\deg P'}\prod_{i=1}^k\big(1-\frac{1}
{q^{\deg \pi_i}}\big)}{q^2\,(q-1)}
\;q^{-2n}\!\!\!\!\sum_{v=(P,Q)\in R^2_{\prim}\;:\;\deg P<\deg Q= q^n,\;P'\,|\,P} 
\Delta_{\frac{(-1)^{n_v}Q_{n_v-1}(v)}{Q_{n_v}(v)}}\;\;\weakstar\;\;
\mu_{Y^{-1}\OOO}\;.
$$
Furthermore, there exists $\tau\in\;]0,\frac{1}{8}]$ such that for
all $\epsilon,\delta>0$, there is a mutiplicative error term in the
above equidistribution claim of the form
$1+\bigO_{\delta,\,P'}(q^{2n(-\tau+\delta)}\,\|g\|_\epsilon)$ when
evaluated on $g$ for every $\epsilon$-locally constant map
$g:Y^{-1}\OOO\ra\RR$.  \ecoro

\dem The result follows by applying the joint equidistribution Theorem
\ref{theo:main} with ${\bf C}=\PP^1$, $\omega=\omega_\infty$ and
$I=P'R$ to the characteristic function of the set
$$
S^1_\infty-S^{1,\sharp}_\infty=\{(x,y)\in\wh K^{\,2}:|x|<|y|=1\}
$$ 
on the right factor, using the following remarks.

\smallskip
$\bullet$~ Let $v=(a,b)\in R^2_{\prim}$ be such that $|a|<|b|$, and
let $(P_i/Q_i)_{-1\leq i\leq n}$ be the sequence of convergents of
$a/b$. Lemma \ref{lem:shortestvsscf} (actually Equation
\eqref{eq:PQcoprime} is sufficient) says that we may take
$w_v=(-(-1)^n\lambda P_{n-1},-(-1)^n\lambda Q_{n-1})$ for $\lambda\in
\FF_q^\times$ such that $v=(\lambda P_n,\lambda Q_n)$. Since
$|P_i|<|Q_i|$ for $0\leq i\leq n$, we have
$$
\frac{z_{w_v}}{z_v}=\frac{y_{w_v}}{y_{v}}=\frac{-(-1)^n\lambda
Q_{n-1}}{\lambda Q_n}=\frac{-(-1)^nQ_{n-1}}{Q_n}\;.
$$

$\bullet$~The map from $Y^{-1}\OOO$ to $\wh K/R$ defined by $f\mapsto
f+R$ is a homeomorphism and an isomorphism of additive groups, which
maps the probability measure $\mu_{Y^{-1}\OOO}$ to $q\;\mu_{\wh K/R}$,
since $\mu_{\wh K/R}$ has total mass $\frac{1}{q}$ by Equation
\eqref{eq:haarKmodR}.

\smallskip
$\bullet$~ The map from $Y^{-1}\OOO$ to itself defined by $f\mapsto
-f$ is an homeomorphism preserving $\mu_{Y^{-1}\OOO}$.

\smallskip
$\bullet$~ We have $\mu_{S^1_\infty}(S^1_\infty-S^{1,\sharp}_\infty)=
\mu_{\wh K}\otimes\mu_{\wh K} (Y^{-1}\OOO\times
\OOO^\times)=\frac{1}{q}\,(1-\frac{1}{q})=\frac{q-1}{q^2}$.

\cqfd

{\small \bibliography{../biblio} }

\bigskip
{\small
\noindent \begin{tabular}{l} 
\\ 
IST Austria, Am Campus 1, 3400 Klosterneuburg, AUSTRA.\\
{\it e-mail: tal.horesh@ist.ac.at}
\end{tabular}
\medskip

\noindent \begin{tabular}{l}
Laboratoire de math\'ematique d'Orsay, UMR 8628 Univ.~Paris-Sud, CNRS\\
Universit\'e Paris-Saclay,
91405 ORSAY Cedex, FRANCE\\
{\it e-mail: frederic.paulin@math.u-psud.fr}
\end{tabular}
}

\end{document}